\documentclass[11pt]{article}%
\usepackage{amsmath}
\usepackage{graphicx}%
\usepackage{amsfonts}%
\usepackage{amssymb}
\newcommand{\eqnb}{\begin{equation}}
\newcommand{\eqne}{\end{equation}}

\begin{document}

\title{\textbf{Tail Probabilities in Queueing Processes}}
\author{Quan-Lin Li\\School of Economics and Management Sciences \\Yanshan University, Qinhuangdao 066004, China }
\maketitle

\begin{abstract}
In the study of large scale stochastic networks with resource management,
differential equations and mean-field limits are two key techniques. Recent
research shows that the expected fraction vector (that is, the tail
probability vector) plays a key role in setting up mean-field differential
equations. To further apply the technique of tail probability vector to deal
with resource management of large scale stochastic networks, this paper
discusses tail probabilities in some basic queueing processes including QBD
processes, Markov chains of GI/M/1 type and of M/G/1 type, and also provides
some effective and efficient algorithms for computing the tail probabilities
by means of the matrix-geometric solution, the matrix-iterative solution, the
matrix-product solution and the two types of $RG$-factorizations. Furthermore,
we consider four queueing examples: The M/M/1 retrial queue, the M(n)/M(n)/1
queue, the M/M/1 queue with server multiple vacations and the M/M/1 queue with
repairable server, where the M/M/1 retrial queue is given a detailed
discussion, while the other three examples are analyzed in less detail. Note
that the results given in this paper will be very useful in the study of large
scale stochastic networks with resource management, including the supermarket
models and the work stealing models.

\vskip    0.5cm


\end{abstract}

\section{Introduction}

We consider a discrete-time (resp. continuous-time) Markov chain whose
transition probability matrix (resp. infinitesimal generator) is given by%
\[
P=\left(
\begin{array}
[c]{ccccc}%
P_{0,0} & P_{0,1} & P_{0,2} & P_{0,3} & \cdots\\
P_{1,0} & P_{1,1} & P_{1,2} & P_{1,3} & \cdots\\
P_{2,0} & P_{2,1} & P_{2,2} & P_{2,3} & \cdots\\
P_{3,0} & P_{3,1} & P_{3,2} & P_{3,3} & \cdots\\
\vdots & \vdots & \vdots & \vdots &
\end{array}
\right)  ,
\]
where the size of the matrix $P_{0,0}$ is $m_{0}$, the size of the matrix
$P_{j,j}$ is $m$ for $j\geq1$, and the sizes of other matrices can be
determined accordingly. We assume that the Markov chain $P$ is irreducible,
aperiodic and positive recurrent. Let $x=\left(  x_{0},x_{1},x_{2}%
,x_{3},\ldots\right)  $ be the stationary probability vector of the Markov
chain $P$, where the size of the vector $x_{0}$ is $m_{0}$ while the size of
the vector $x_{j}$ is $m$ for $j\geq1$. The main purpose of this paper is to
discuss the tail probabilities: $\pi_{k}=\sum_{j=k}^{\infty}x_{j}$ and to
provide some efficient algorithms for computing the tail probabilities
$\pi_{k}$ for $k\geq1$.

Recent queueing literature indicates that the study of tail probabilities
$\left\{  \pi_{k},k\geq0\right\}  $ plays a key role in analyzing large scale
stochastic networks with resource management, such as, the supermarket models
and the work stealing models, e.g., see Vvedenskaya and Suhov \cite{Vve:1997}
and Mitzenmacher \cite{Mit:1999}. When considering a large scale stochastic
network with resource management, differential equations and mean-field limits
are always two key techniques, while the tail probabilities play a key role in
setting up mean-field differential equations. The detailed interpretation on
the mean-field differential equations was given in Vvedenskaya et al
\cite{Vve:1996}, Mitzenmacher \cite{Mit:1996}, Ethier and Kurtz
\cite{Eth:1986} and Kurtz \cite{Kur:1981}. In the first two papers, the
authors considered a supermarket model with $N$ identical servers, where the
service times are exponential with service rate $\mu$, and the input flow is
Poisson with arrival rate $N\lambda$. Upon arrival, each customer chooses
$d\geq1$ servers from the $N$ servers independently and uniformly at random,
and joins the one whose queue length is the shortest. Let $n_{k}^{(N)}(t)$
denote the number of servers queued by at least $k\geq0$ customers at time
$t$, and $u_{k}\left(  t\right)  =\lim_{N\rightarrow\infty}E\left[
n_{k}^{(N)}(t)/N\right]  $. If $\rho=\lambda/\mu<1$, then the supermarket
model is stable, and%
\begin{equation}
\frac{\text{d}}{\text{d}t}u_{k}\left(  t\right)  =\lambda\left\{  \left[
u_{k-1}\left(  t\right)  \right]  ^{d}-\left[  u_{k}\left(  t\right)  \right]
^{d}\right\}  -\mu\left[  u_{k}\left(  t\right)  -u_{k+1}\left(  t\right)
\right]  \label{IN-1}%
\end{equation}
with the boundary condition $u_{0}\left(  t\right)  =1$. We write that
$\pi_{k}=\lim_{t\rightarrow+\infty}u_{k}\left(  t\right)  $ for $k\geq0$. Then
$\pi_{0}=1$ and for $k\geq1$%
\begin{equation}
\lambda\left(  \pi_{k-1}^{d}-\pi_{k}^{d}\right)  -\mu\left(  \pi_{k}-\pi
_{k+1}\right)  =0.\label{IN-2}%
\end{equation}
This gives%
\[
\pi_{k}=\rho^{\frac{d^{k}-1}{d-1}},\text{ \ \ }k\geq1.
\]
Specifically, $\pi_{1}=\rho$ is directly derived by%
\[
\lambda\sum_{k=1}^{\infty}\left(  \pi_{k-1}^{d}-\pi_{k}^{d}\right)  -\mu
\sum_{k=1}^{\infty}\left(  \pi_{k}-\pi_{k+1}\right)  =0.
\]
If $d=1$, then $\pi_{k}=\rho^{k}$ for $k\geq0$ are the tail probabilities of
the M/M/1 queue.

Since the introduction of the expected fraction vector (or the tail
probability vector) by Vvedenskaya et al \cite{Vve:1996} and Mitzenmacher
\cite{Mit:1996}, research on supermarket models and work stealing models has
been greatly motivated by some practical applications such as computer
networks, manufacturing systems and transportation networks. Subsequent papers
have been published on this theme, among which, see, modeling more crucial
factors by Mitzenmacher \cite{Mit:1998, Mit:1999}, Jacquet and Vvedenskaya
\cite{Jac:1998}, Jacquet et al \cite{Jac:1999} and Vvedenskaya and Suhov
\cite{Vve:2005}; studying fast Jackson networks by\ Martin and Suhov
\cite{Mar:1999}, Martin \cite{Mar:2001} and Suhov and Vvedenskaya
\cite{Suh:2002}; discussing value of information by Mitzenmacher
\cite{Mit:2000} and Mitzenmacher et al \cite{Mit:2002}; analyzing
non-exponential server times and/or non-Poisson inputs by Mitzenmacher
\cite{Mit:1996}, Vvedenskaya and Suhov \cite{Vve:1997}, Bramson
\cite{Bra:2011a}, Bramson et al \cite{Bra:2010, Bra:2011, Bra:2012}, Li et al
\cite{LiLW:2011}, Li and Lui \cite{LiL:2010} and Li \cite{Li:2011}. For a
comprehensive analysis of supermarket models and work stealing models, readers
may refer to Vvedenskaya and Suhov \cite{Vve:1997}, Mitzenmacher et al
\cite{Mit:2001} and Mitzenmacher and Upfal \cite{Mit:2005}. From those papers,
it is seen that the tail probabilities $\left\{  \pi_{k},k\geq0\right\}  $ is
obtained from the mean-field differential equations as $N\rightarrow\infty$
and $t\rightarrow+\infty$, and also it is a key to analyze performance
measures of the supermarket models and of the work stealing models.

During the last two decades considerable attention has been paid to studying
QBD processes, which has been well documented, for example, by Chapter 3 of
Neuts \cite{Neu:1981}, Naoumov \cite{Nao:1996}, Bright and Taylor
\cite{Bri:1995, Bri:1996}, Ramaswami \cite{Ram:1996}, Latouche and Ramaswami
\cite{Lat:1999} and Li and Cao \cite{Li:2004}. For Markov chains of GI/M/1
type and Markov chains of M/G/1 type, readers may refer to four excellent
books by Neuts \cite{Neu:1981, Neu:1989}, Latouche and Ramaswami
\cite{Lat:1999} and Li \cite{Li:2010}.

Some papers were published on asymptotic behavior of the stationary
probability vectors for both queueing systems and Markov chains. Readers may
refer to, such as, Markov chains of $GI/M/1$ type by Neuts \cite{Neu:1981};
Markov chains of $M/G/1$ type by Falkenberg \cite{Fal:1994}, Abate et al
\cite{Aba:1994a}, Choudhury and Whitt \cite{Chou:1994}, Asmussen and M\o ller
\cite{Asm:1999} and Takine \cite{Tak:2004}; and Markov chains of $GI/G/1$ type
by Li and Zhao \cite{Li:2005a, Li:2005b}.

The main purpose of this paper is to provide some novel and efficient
algorithms for computing the tail probabilities in three classes of important
Markov chains: QBD processes, Markov chains of GI/M/1 type and Markov chains
of M/G/1 type. Note that the algorithms are based on the matrix-geometric
solution, the matrix-iterative solution, the matrix-product solution and the
two types of $RG$-factorizations. Also, we consider four queueing examples:
The M/M/1 retrial queue, the M(n)/M(n)/1 queue, the M/M/1 queue with server
multiple vacations and the M/M/1 queue with repairable server. Based on this,
it is seen that the method of this paper can deal with more general queue
examples, such as, the MAP/PH/1 queue, the GI/PH/1 queue and the BMAP/SM/1
queue. Therefore, the results of this paper are very useful in setting up the
mean-field differential equations for large scale stochastic networks with
resource management, including the supermarket models and the work stealing models.

The remainder of this paper is organized as follows. In Section 2, we analyze
a continuous-time level-independent QBD process. When the QBD process is
irreducible, aperiodic and positive recurrent, we apply the matrix-geometric
solution and the two types of $RG$-factorizations to compute the tail
probabilities in the stationary regime. In Section 3, we consider an
continuous-time level-dependent QBD process, and compute the tail
probabilities in the stationary regime. In Section 4, we discuss two classes
of important Markov chains: Markov chains of GI/M/1 type and of M/G/1 type,
and derive the tail probabilities in the stationary regime. In Section 5, we
study four queueing examples, where the M/M/1 retrial queue is given a
detailed discussion, while the other three queues are analyzed in less detail.
Some concluding remarks are given in Section 6.

\section{Level-Independent QBD Processes}

In this section, we consider a continuous-time level-independent QBD process.
When the QBD process is irreducible, aperiodic and positive recurrent, we
apply the two types of $RG$-factorizations to compute the tail probabilities
in the stationary regime. Furthermore, the tail probabilities of the
stationary probability vector is well related to the matrix-geometric solution
when the UL-type $RG$-factorization is used.

We consider a continuous-time level-independent QBD process whose
infinitesimal generator is given by
\begin{equation}
Q=\left(
\begin{array}
[c]{cccccc}%
B_{1} & B_{0} &  &  &  & \\
B_{2} & A_{1} & A_{0} &  &  & \\
& A_{2} & A_{1} & A_{0} &  & \\
&  & A_{2} & A_{1} & A_{0} & \\
&  &  & \ddots & \ddots & \ddots
\end{array}
\right)  , \label{LI-1}%
\end{equation}
where the sizes of the two matrices $B_{1}$ and $A_{1}$ are $m_{0}$ and $m$,
respectively; and the sizes of other matrices can be determined accordingly.
We assume that this QBD process is irreducible, aperiodic and positive recurrent.

\subsection{The matrix-geometric solution}

Let $x=\left(  x_{0},x_{1},x_{2},\ldots\right)  $ be the stationary
probability vector of the QBD process, and $R$ and $G$ the minimal nonnegative
solutions to the nonlinear equations $A_{0}+RA_{1}+R^{2}A_{2}=0$ and
$A_{0}G^{2}+A_{1}G+A_{2}=0$, respectively. Then%
\begin{equation}
x_{k}=x_{1}R^{k-1},\text{ \ \ }k\geq2, \label{LI-2}%
\end{equation}
where $x_{0}$ and $x_{1}$ are uniquely determined by the following system of
linear equations%
\[
\left\{
\begin{array}
[c]{l}%
x_{0}B_{1}+x_{1}B_{2}=0,\\
x_{0}B_{0}+x_{1}\left(  A_{1}+RA_{2}\right)  =0,\\
x_{0}e+x_{1}\left(  I-R\right)  ^{-1}e=1,
\end{array}
\right.
\]
where $e$ is a column vector of ones.

We write%
\begin{equation}
\pi_{k}=\sum_{j=k}^{\infty}x_{j},\text{ \ \ }k\geq1. \label{LI-3}%
\end{equation}
It follows from (\ref{LI-2}) that%
\begin{equation}
\pi_{k}=x_{1}\left(  I-R\right)  ^{-1}R^{k-1},\text{ \ \ }k\geq1. \label{LI-4}%
\end{equation}

\subsection{The UL-type $RG$-factorization}

Now, we apply the UL-type $RG$-factorization to provide a novel method for
deriving the tail probability vector $\pi=\left(  \pi_{1},\pi_{2},\pi
_{3},\ldots\right)  $.

Note that $xQ=0$, so we have%
\begin{equation}
\left\{
\begin{array}
[c]{ll}%
x_{0}B_{0}+x_{1}A_{1}+x_{2}A_{2}=0, & k=1,\\
x_{k-1}A_{0}+x_{k}A_{1}+x_{k+1}A_{2}=0, & k\geq2.
\end{array}
\right.  \label{LI-5}%
\end{equation}
This gives%
\[
\left\{
\begin{array}
[c]{ll}%
\pi_{1}\left(  A_{0}+A_{1}\right)  +\pi_{2}A_{2}=-x_{0}B_{0}, & k=1,\\
\pi_{k-1}A_{0}+\pi_{k}A_{1}+\pi_{k+1}A_{2}=0, & k\geq2.
\end{array}
\right.
\]
Hence we obtain%
\begin{equation}
\pi\mathbb{Q}=-\left(  x_{0}B_{0},0,0,0,\ldots\right)  , \label{LI-6}%
\end{equation}
where%
\[
\mathbb{Q=}\left(
\begin{array}
[c]{cccccc}%
A_{0}+A_{1} & A_{0} &  &  &  & \\
A_{2} & A_{1} & A_{0} &  &  & \\
& A_{2} & A_{1} & A_{0} &  & \\
&  & A_{2} & A_{1} & A_{0} & \\
&  &  & \ddots & \ddots & \ddots
\end{array}
\right)  .
\]

Let%
\[
\Phi_{0}=\left(  A_{0}+A_{1}\right)  +RA_{2},
\]%
\[
\Phi_{k}=\Phi=A_{1}+RA_{2},\text{ \ \ }k\geq1.
\]
Then the UL-type $RG$-factorization of the matrix $\mathbb{Q}$ is given by%
\begin{equation}
\mathbb{Q=}\left(  I-R_{U}\right)  U\left(  I-G_{L}\right)  , \label{LI-7}%
\end{equation}
where%
\[
R_{U}=\left(
\begin{array}
[c]{ccccc}%
0 & R &  &  & \\
& 0 & R &  & \\
&  & 0 & R & \\
&  &  & \ddots & \ddots
\end{array}
\right)  ,
\]%
\[
U=\text{diag}\left(  \Phi_{0},\Phi,\Phi,\ldots\right)
\]
and%
\[
G_{L}=\left(
\begin{array}
[c]{ccccc}%
0 &  &  &  & \\
G & 0 &  &  & \\
& G & 0 &  & \\
&  & G & 0 & \\
&  &  & \ddots & \ddots
\end{array}
\right)  .
\]

It follows from (\ref{LI-6}) and (\ref{LI-7}) that%
\begin{align*}
\pi &  =-\left(  x_{0}B_{0},0,0,0,\ldots\right)  \left(  I-G_{L}\right)
^{-1}U^{-1}\left(  I-R_{U}\right)  ^{-1}\\
&  =-\left(  x_{0}B_{0}\Phi_{0}^{-1},0,0,0,\ldots\right)  \left(
I-R_{U}\right)  ^{-1}.
\end{align*}
Note that%
\[
\left(  I-R_{U}\right)  ^{-1}=\left(
\begin{array}
[c]{ccccc}%
I & R & R^{2} & R^{3} & \cdots\\
& I & R & R^{2} & \cdots\\
&  & I & R & \cdots\\
&  &  & I & \cdots\\
&  &  &  & \ddots
\end{array}
\right)  ,
\]
and so we obtain%
\begin{equation}
\left\{
\begin{array}
[c]{ll}%
\pi_{1}=x_{0}B_{0}\left(  -\Phi_{0}^{-1}\right)  , & k=1,\\
\pi_{k}=x_{0}B_{0}\left(  -\Phi_{0}^{-1}\right)  R^{k-1}, & k\geq2.
\end{array}
\right.  \label{LI-8}%
\end{equation}

Comparing (\ref{LI-8}) with (\ref{LI-4}), we obtain%
\begin{equation}
x_{0}B_{0}\left(  -\Phi_{0}^{-1}\right)  =x_{1}\left(  I-R\right)  ^{-1}.
\label{LI-9}%
\end{equation}
This gives%
\[
x_{1}=x_{0}B_{0}\left(  -\Phi_{0}^{-1}\right)  \left(  I-R\right)  .
\]

\subsection{The LU-type $RG$-factorization}

In what follows we apply the LU-type $RG$-factorization to provide a novel and
effective method for deriving the tail probability vector $\pi=\left(  \pi
_{1},\pi_{2},\pi_{3},\ldots\right)  $.

Let%
\begin{equation}
\Psi_{0}=A_{0}+A_{1} \label{LU-1}%
\end{equation}
and for $k\geq1$%
\begin{equation}
\Psi_{k}=A_{1}+A_{2}\left(  -\Psi_{k-1}^{-1}\right)  A_{0}. \label{LU-2}%
\end{equation}
We write that for $k\geq1$%
\begin{equation}
\mathbf{R}_{k}=A_{2}\left(  -\Psi_{k-1}^{-1}\right)  \label{LU-3}%
\end{equation}
and%
\begin{equation}
\mathbf{G}_{k-1}=\left(  -\Psi_{k-1}^{-1}\right)  A_{0}. \label{LU-4}%
\end{equation}
Then the LU-type $RG$-factorization of the matrix $\mathbb{Q}$ is given by%
\begin{equation}
\mathbb{Q=}\left(  I-R_{L}\right)  \mathbf{U}\left(  I-G_{U}\right)  ,
\label{LI-10}%
\end{equation}
where%
\[
R_{L}=\left(
\begin{array}
[c]{ccccc}%
0 &  &  &  & \\
\mathbf{R}_{1} & 0 &  &  & \\
& \mathbf{R}_{2} & 0 &  & \\
&  & \mathbf{R}_{3} & 0 & \\
&  &  & \ddots & \ddots
\end{array}
\right)  ,
\]%
\[
U=\text{diag}\left(  \Psi_{0},\Psi_{1},\Psi_{2},\ldots\right)
\]
and%
\[
G_{U}=\left(
\begin{array}
[c]{ccccc}%
0 & \mathbf{G}_{0} &  &  & \\
& 0 & \mathbf{G}_{1} &  & \\
&  & 0 & \mathbf{G}_{2} & \\
&  &  & \ddots & \ddots
\end{array}
\right)  .
\]
Let%
\[
X_{k}^{\left(  l\right)  }=\mathbf{R}_{l}\mathbf{R}_{l-1}\mathbf{R}%
_{l-2}\cdots\mathbf{R}_{l-k+1},\text{ \ }l\geq k\geq1,
\]%
\[
Y_{k}^{\left(  l\right)  }=\mathbf{G}_{l}\mathbf{G}_{l+1}\mathbf{G}%
_{l+2}\cdots\mathbf{G}_{l+k-1,\text{ \ \ }}k\geq1,l\geq0.
\]
Then%
\[
-\mathbf{U}^{-1}=\text{diag}\left(  -\Psi_{0}^{-1},-\Psi_{1}^{-1},-\Psi
_{2}^{-1},-\Psi_{3}^{-1},\ldots\right)  ,
\]%
\[
\left(  I-R_{L}\right)  ^{-1}=\left(
\begin{array}
[c]{ccccc}%
I &  &  &  & \\
X_{1}^{\left(  1\right)  } & I &  &  & \\
X_{2}^{\left(  2\right)  } & X_{1}^{\left(  2\right)  } & I &  & \\
X_{3}^{\left(  3\right)  } & X_{2}^{\left(  3\right)  } & X_{1}^{\left(
3\right)  } & I & \\
\vdots & \vdots & \vdots & \vdots & \ddots
\end{array}
\right)  ,
\]%
\[
\left(  I-G_{U}\right)  ^{-1}=\left(
\begin{array}
[c]{ccccc}%
I & Y_{1}^{\left(  0\right)  } & Y_{2}^{\left(  0\right)  } & Y_{3}^{\left(
0\right)  } & \cdots\\
& I & Y_{1}^{\left(  1\right)  } & Y_{2}^{\left(  1\right)  } & \cdots\\
&  & I & Y_{1}^{\left(  2\right)  } & \cdots\\
&  &  & I & \cdots\\
&  &  &  & \ddots
\end{array}
\right)  .
\]
It follows from (\ref{LI-6}) and (\ref{LI-7}) that%
\[
\pi=-\left(  x_{0}B_{0},0,0,0,\ldots\right)  \left(  I-G_{U}\right)
^{-1}\mathbf{U}^{-1}\left(  I-R_{L}\right)  ^{-1},
\]
this gives%
\begin{equation}
\pi_{1}=x_{0}B_{0}\left[  \left(  -\Psi_{0}^{-1}\right)  +\sum_{k=1}^{\infty
}Y_{k}^{\left(  0\right)  }\left(  -\Psi_{k}^{-1}\right)  X_{k}^{\left(
k\right)  }\right]  , \label{LI-11}%
\end{equation}
and $n\geq2$%
\begin{equation}
\pi_{n}=x_{0}B_{0}\left[  Y_{n-1}^{\left(  0\right)  }\left(  -\Psi_{n-1}%
^{-1}\right)  +\sum_{k=n}^{\infty}Y_{k}^{\left(  0\right)  }\left(  -\Psi
_{k}^{-1}\right)  X_{k-\left(  n-1\right)  }^{\left(  k\right)  }\right]  .
\label{LI-12}%
\end{equation}

The expressions (\ref{LI-11}) and (\ref{LI-12}) for the tail probability
vector $\left\{  \pi_{k}:k\geq1\right\}  $ seem complicated, but they can
easily be computed by means of the iterative relations (\ref{LU-1}) to
(\ref{LU-4}) through some simple matrix calculations.

\section{Level-Dependent QBD Processes}

In this section, we consider a continuous-time level-dependent QBD process.
When the QBD process is irreducible, aperiodic and positive recurrent, we
apply the LU-type $RG$-factorization to compute the tail probabilities in the
stationary regime. Similarly, we can apply the UL-type $RG$-factorization to
compute the tail probabilities without any difficulty.

We consider a continuous-time level-dependent QBD process whose infinitesimal
generator is given by
\begin{equation}
Q=\left(
\begin{array}
[c]{ccccc}%
A_{1}^{(0)} & A_{0}^{(0)} &  &  & \\
A_{2}^{(1)} & A_{1}^{(1)} & A_{0}^{(1)} &  & \\
& A_{2}^{(2)} & A_{1}^{(2)} & A_{0}^{(2)} & \\
&  & \ddots & \ddots & \ddots
\end{array}
\right)  , \label{LD-1}%
\end{equation}
where the size of the matrix $A_{1}^{(0)}$ is $m_{0}$ while the size of the
matrix $A_{1}^{(k)}$ is $m$ for $k\geq1$. We assume that this QBD process is
irreducible, aperiodic and positive recurrent.

\subsection{The matrix-product solution}

Let the matrix sequence $\left\{  R_{l},l\geq0\right\}  $ be the minimal
nonnegative solution to the system of nonlinear matrix equations%
\[
A_{0}^{(l)}+R_{l}A_{1}^{(l+1)}+R_{l}R_{l+1}A_{2}^{(l+2)}=0,\ \ \ l\geq0.
\]
Using Chapter 1 of Li \cite{Li:2010}, we have%
\begin{equation}
x_{0}=\kappa v, \label{LD-6}%
\end{equation}%
\begin{equation}
x_{k}=\kappa vR_{0}R_{1}\cdots R_{k-1},\ \ k\geq1, \label{LD-7}%
\end{equation}
where $v$ is the stationary probability vector of the censored chain
$U_{0}=A_{1}^{\left(  0\right)  }+R_{0}A_{2}^{\left(  1\right)  }$ to level
$0$, and the constant $\kappa$ is given by%
\[
\kappa=\frac{1}{1+v\left(  \sum\limits_{k=0}^{\infty}R_{0}R_{1}\cdots
R_{k}\right)  e}%
\]
Therefore, it follows from (\ref{LD-6}) and (\ref{LD-7}) that%
\[
\pi_{k}=\sum_{j=k}^{\infty}x_{j}=\kappa vR_{0}R_{1}\cdots R_{k-1}\left(
I+\sum_{l=0}^{\infty}R_{k}R_{k+1}\cdots R_{k+l}\right)  .
\]

\subsection{The LU-type $RG$-factorization}

Here, we only provide a detailed analysis for applying the LU-type
$RG$-factorization to compute the tail probabilities, while the UL-type
$RG$-factorization can be used similarly for such an analysis.

Since%
\[
\pi_{k}=\sum_{j=k}^{\infty}x_{j},\text{ \ \ }k\geq1,
\]
it is easy to see that $x_{k}=\pi_{k}-\pi_{k+1}$ for $k\geq1$. Note that
$xQ=0$, we have%
\[
\left\{
\begin{array}
[c]{ll}%
\pi_{1}A_{1}^{(1)}+\pi_{2}\left[  A_{2}^{(2)}-A_{1}^{(1)}\right]  -\pi
_{3}A_{2}^{(2)}=-x_{0}A_{0}^{(0)}, & k=1,\\
\pi_{k-1}A_{0}^{(k-1)}+\pi_{k}\left[  A_{1}^{(k)}-A_{0}^{(k-1)}\right]
+\pi_{k+1}\left[  A_{2}^{(k+1)}-A_{1}^{(k)}\right]  -\pi_{k+2}A_{2}%
^{(k+1)}=0, & k\geq2,
\end{array}
\right.
\]
we obtain that $\pi\mathbf{Q}=-\left(  x_{0}A_{0}^{(0)},0,0,0,\ldots\right)
$, where%
\[
\mathbf{Q=}\left(
\begin{array}
[c]{ccccccc}%
A_{1}^{(1)} & A_{0}^{(1)} &  &  &  &  & \\
A_{2}^{(2)}-A_{1}^{(1)} & A_{1}^{(2)}-A_{0}^{(1)} & A_{0}^{(2)} &  &  &  & \\
-A_{2}^{(2)} & A_{2}^{(3)}-A_{1}^{(2)} & A_{1}^{(3)}-A_{0}^{(2)} & A_{0}^{(3)}
&  &  & \\
& -A_{2}^{(3)} & A_{2}^{(4)}-A_{1}^{(3)} & A_{1}^{(4)}-A_{0}^{(3)} &
A_{0}^{(4)} &  & \\
&  & -A_{2}^{(4)} & A_{2}^{(5)}-A_{1}^{(4)} & A_{1}^{(5)}-A_{0}^{(4)} &
A_{0}^{(5)} & \\
&  &  & \ddots & \ddots & \ddots & \ddots
\end{array}
\right)  .
\]
This gives%
\begin{equation}
\pi\left(  \widehat{Q}-\mathbb{Q}\right)  =-\left(  x_{0}A_{0}^{(0)}%
,0,0,0,\ldots\right)  , \label{LD-3}%
\end{equation}
where%
\[
\widehat{Q}=\left(
\begin{array}
[c]{ccccc}%
A_{1}^{(1)} & A_{0}^{(1)} &  &  & \\
A_{2}^{(2)} & A_{1}^{(2)} & A_{0}^{(2)} &  & \\
& A_{2}^{(3)} & A_{1}^{(3)} & A_{0}^{(3)} & \\
&  & \ddots & \ddots & \ddots
\end{array}
\right)
\]
and $\mathbf{0}=\left(  0,0,0,0,\ldots\right)  $, $0$ is an $m\times m$ zero
matrix,%
\[
\mathbb{Q=}\left(
\begin{array}
[c]{c}%
\mathbf{0}\\
\widehat{Q}%
\end{array}
\right)  .
\]
It follows from (\ref{LD-3}) that%
\[
\pi\left(  I-\mathbb{Q}\widehat{Q}_{\max}^{-1}\right)  =-\left(  x_{0}%
A_{0}^{(0)},0,0,0,\ldots\right)  \widehat{Q}_{\max}^{-1},
\]
where $\widehat{Q}_{\max}^{-1}$ is the maximal non-positive inverse of the
infinitesimal generator $\widehat{Q}$. Hence, this gives%
\[
\pi=-\left(  x_{0}A_{0}^{(0)},0,0,0,\ldots\right)  \widehat{Q}_{\max}^{-1}%
\sum_{k=0}^{\infty}\left(  \mathbb{Q}\widehat{Q}_{\max}^{-1}\right)  ^{k}.
\]

Now, we apply the LU-type $RG$-factorization to provide the maximal
non-positive inverse $\widehat{Q}_{\max}^{-1}$ of the infinitesimal generator
$\widehat{Q}$. To that end, we write%
\[
\Psi_{0}=A_{1}^{(1)}%
\]
and for $k\geq1$%
\[
\Psi_{k}=A_{1}^{(k+1)}+A_{2}^{\left(  k+1\right)  }\left(  -\Psi_{k-1}%
^{-1}\right)  A_{0}^{(k)}.
\]
It is easy to check that $\Psi_{l}$ is the infinitesimal generator of an
irreducible continuous-time Markov chain, and the Markov chain $\Psi_{l}$ is
transient. Thus the matrix $\Psi_{l}$ is invertible for $l\geq0$.

Based on the $U$-measure $\left\{  \Psi_{l}\right\}  $, for $k\geq1$ we can
respectively define the LU-type $R$- and $G$-measures as%
\[
\mathbf{R}_{k}=A_{2}^{\left(  k+1\right)  }\left(  -\Psi_{k-1}^{-1}\right)
\]
and%
\[
\mathbf{G}_{k-1}=\left(  -\Psi_{k-1}^{-1}\right)  A_{0}^{(k)}.
\]

Note that the matrix sequence $\left\{  \mathbf{R}_{k}:k\geq1\right\}  $ is
the unique nonnegative solution to the system of nonlinear matrix equations%
\[
\mathbf{R}_{k+1}\mathbf{R}_{k}A_{0}^{\left(  k\right)  }+\mathbf{R}_{k+1}%
A_{1}^{\left(  k+1\right)  }+A_{2}^{\left(  k+2\right)  }=0,
\]
with the boundary condition%
\[
\mathbf{R}_{1}=A_{2}^{\left(  2\right)  }\left(  -\Psi_{0}^{-1}\right)  .
\]
Hence we obtain%
\[
\mathbf{R}_{k+1}=-A_{2}^{\left(  k+2\right)  }\left[  \mathbf{R}_{k}%
A_{0}^{\left(  k\right)  }+A_{1}^{\left(  k+1\right)  }\right]  ^{-1}.
\]
Similarly, the matrix sequence $\{\mathbf{G}_{k}:k\geq0\}$ is the unique
nonnegative solution to the system of nonlinear matrix equations%
\[
A_{0}^{\left(  k+1\right)  }+A_{1}^{\left(  k+1\right)  }\mathbf{G}_{k}%
+A_{2}^{\left(  k+1\right)  }\mathbf{G}_{k-1}\mathbf{G}_{k}=0,
\]
with the boundary condition%
\[
\mathbf{G}_{0}=\left(  -\Psi_{0}^{-1}\right)  A_{0}^{(1)}.
\]
Thus we obtain%
\[
\mathbf{G}_{k}=-\left[  A_{1}^{\left(  k+1\right)  }+A_{2}^{\left(
k+1\right)  }\mathbf{G}_{k-1}\right]  ^{-1}A_{0}^{(k+1)}.
\]
The LU-type $RG$-factorization of the QBD process $\widehat{Q}$ is given by%
\begin{equation}
\widehat{Q}=\left(  I-\mathbb{R}_{L}\right)  \mathbb{U}_{D}\left(
I-\mathbb{G}_{U}\right)  , \label{LD-2}%
\end{equation}
where%
\[
\mathbb{R}_{L}=\left(
\begin{array}
[c]{ccccc}%
0 &  &  &  & \\
\mathbf{R}_{1} & 0 &  &  & \\
& \mathbf{R}_{2} & 0 &  & \\
&  & \mathbf{R}_{3} & 0 & \\
&  &  & \ddots & \ddots
\end{array}
\right)  ,
\]%
\[
\mathbb{U}_{D}=\text{diag}\left(  \Psi_{0},\Psi_{1},\Psi_{2},\Psi_{3}%
,\ldots\right)  ,
\]%
\[
\mathbb{G}_{U}=\left(
\begin{array}
[c]{ccccc}%
0 & \mathbf{G}_{0} &  &  & \\
& 0 & \mathbf{G}_{1} &  & \\
&  & 0 & \mathbf{G}_{2} & \\
&  &  & 0 & \ddots\\
&  &  &  & \ddots
\end{array}
\right)  .
\]
Let%
\begin{equation}
X_{k}^{\left(  l\right)  }=\mathbf{R}_{l}\mathbf{R}_{l-1}\mathbf{R}%
_{l-2}\cdots\mathbf{R}_{l-k+1},\text{ \ }l\geq k\geq1, \label{LD-4}%
\end{equation}
and%
\begin{equation}
Y_{k}^{\left(  l\right)  }=\mathbf{G}_{l}\mathbf{G}_{l+1}\mathbf{G}%
_{l+2}\cdots\mathbf{G}_{l+k-1},\text{ \ }k\geq1,l\geq0. \label{LD-5}%
\end{equation}
Then%
\[
\left(  I-\mathbb{R}_{L}\right)  ^{-1}=\left(
\begin{array}
[c]{ccccc}%
I &  &  &  & \\
X_{1}^{\left(  1\right)  } & I &  &  & \\
X_{2}^{\left(  2\right)  } & X_{1}^{\left(  2\right)  } & I &  & \\
X_{3}^{\left(  3\right)  } & X_{2}^{\left(  3\right)  } & X_{1}^{\left(
3\right)  } & I & \\
\vdots & \vdots & \vdots & \vdots & \ddots
\end{array}
\right)
\]
and%
\[
\left(  I-\mathbb{G}_{U}\right)  ^{-1}=\left(
\begin{array}
[c]{ccccc}%
I & Y_{1}^{\left(  0\right)  } & Y_{2}^{\left(  0\right)  } & Y_{3}^{\left(
0\right)  } & \cdots\\
& I & Y_{1}^{\left(  1\right)  } & Y_{2}^{\left(  1\right)  } & \cdots\\
&  & I & Y_{1}^{\left(  2\right)  } & \cdots\\
&  &  & I & \cdots\\
&  &  &  & \ddots
\end{array}
\right)  .
\]
Hence we obtain%
\begin{align*}
\pi &  =-\left(  x_{0}A_{0}^{(1)},0,0,0,\ldots\right)  \widehat{Q}_{\max}%
^{-1}\sum_{k=0}^{\infty}\left(  \mathbb{Q}\widehat{Q}_{\max}^{-1}\right)
^{k}\\
&  =-\left(  x_{0}A_{0}^{(1)},0,0,0,\ldots\right)  \left(  I-\mathbb{G}%
_{U}\right)  ^{-1}\mathbb{U}_{D}^{-1}\left(  I-\mathbb{R}_{L}\right)
^{-1}\sum_{k=0}^{\infty}\left[  \mathbb{Q}\left(  I-\mathbb{G}_{U}\right)
^{-1}\mathbb{U}_{D}^{-1}\left(  I-\mathbb{R}_{L}\right)  ^{-1}\right]  ^{k}.
\end{align*}
This can be calculated by some ordinary matrix computation.

\section{Two Classes of Important Markov Chains}

In this section, we consider two classes of important Markov chains: Markov
chains of GI/M/1 type and of M/G/1 type, each of which is basic in the study
of queueing processes, e.g., see Neuts \cite{Neu:1981, Neu:1989} for more
details. We provide two different methods to derive the tail probabilities of
stationary probability vectors of the two classes of Markov chains.

\subsection{Markov chains of GI/M/1 type}

We consider a discrete-time Markov chain $P$ of $GI/M/1$ type whose transition
matrix is given by
\begin{equation}
P=\left(
\begin{array}
[c]{cccccc}%
B_{1} & B_{0} &  &  &  & \\
B_{2} & A_{1} & A_{0} &  &  & \\
B_{3} & A_{2} & A_{1} & A_{0} &  & \\
B_{4} & A_{3} & A_{2} & A_{1} & A_{0} & \\
\vdots & \vdots & \vdots & \vdots & \vdots & \ddots
\end{array}
\right)  , \label{D-1}%
\end{equation}
where the sizes of the two matrices $B_{1}$ and $A_{1}$ are $m_{0}$ and $m$,
respectively, while the sizes of other matrices can be determined accordingly.
We assume that this Markov chain is irreducible, aperiodic and positive
recurrent. Let the matrix $R$ be the minimal nonnegative solution to the
nonlinear matrix equation $R=\sum_{k=0}^{\infty}R^{k}A_{k}$.

In what follows we provide two methods to derive the tail probabilities in the
stationary regime.

\textbf{(a) The matrix-geometric solution}

Using Chapter 2 of Li \cite{Li:2010}, the stationary probability vector
$x=\left(  x_{0},x_{1},x_{2},\ldots\right)  $ is given by%

\[
\left\{
\begin{array}
[c]{l}%
x_{0}=\tau y_{0},\\
x_{k}=x_{0}R_{1}R^{k-1},\ \ k\geq1,
\end{array}
\right.
\]
where%
\[
R_{1}=\left(  I-\sum\limits_{k=0}^{\infty}R^{k}B_{k+1}\right)  ^{-1}B_{0}%
\]
and%
\[
\Psi_{0}=\sum\limits_{k=0}^{\infty}R^{k}B_{k+1},
\]
$y_{0}$ is the stationary probability vector of the censored Markov chain
$\Psi_{0}$ to level $0$,\ and the scalar $\tau$ is determined by%
\[
\tau=\frac{1}{1+y_{0}R_{1}\left(  I-R\right)  ^{-1}e}.
\]
Thus, for $k\geq1$ we have%
\begin{equation}
\pi_{k}=\sum_{j=k}^{\infty}x_{j}=x_{0}R_{1}\left(  I-R\right)  ^{-1}R^{k-1}.
\label{D-2}%
\end{equation}

\textbf{(b) The UL-type }$\mathbf{RG}$\textbf{-factorization}

Note that $x=xP$, we obtain%
\[
\left\{
\begin{array}
[c]{ll}%
\pi_{1}=\pi_{1}\left(  A_{0}+A_{1}\right)  +\sum\limits_{k=2}^{\infty}\pi
_{k}A_{k}+x_{0}B_{0}, & k=1,\\
\pi_{k}=\sum\limits_{j=0}^{\infty}\pi_{k-1+j}A_{j}, & k\geq2,
\end{array}
\right.
\]
this gives%
\begin{equation}
\pi=\pi\mathbb{P}+\left(  x_{0}B_{0},0,0,0,\ldots\right)  , \label{D-3}%
\end{equation}
where%
\[
\mathbb{P=}\left(
\begin{array}
[c]{cccccc}%
A_{0}+A_{1} & A_{0} &  &  &  & \\
A_{2} & A_{1} & A_{0} &  &  & \\
A_{3} & A_{2} & A_{1} & A_{0} &  & \\
A_{4} & A_{3} & A_{2} & A_{1} & A_{0} & \\
\vdots & \vdots & \vdots & \vdots & \vdots & \ddots
\end{array}
\right)  ,
\]
which is of GI/M/1 type. Then using Chapter 2 of Li \cite{Li:2010}, the
$U$-measure is given by%
\[
\widehat{\Psi}_{0}=\left(  A_{0}+A_{1}\right)  +\sum_{k=2}^{\infty}%
R_{1}R^{k-2}A_{k}%
\]
and for $k\geq1$%
\[
\widehat{\Psi}=\widehat{\Psi}_{k}=\sum_{k=1}^{\infty}R^{k-1}A_{k};
\]
the $R$-measure is given by%
\[
R_{k}=R,\text{ \ }k\geq1,
\]
and the $G$-measure%
\[
G_{j,0}=\left(  I-\widehat{\Psi}\right)  ^{-1}\left(  \sum_{k=j+1}^{\infty
}R^{k-1}A_{k}\right)  ,\ \ j\geq1,
\]
and%
\[
G_{j}=\left(  I-\widehat{\Psi}\right)  ^{-1}\left(  \sum_{k=j+1}^{\infty
}R^{k-1}A_{k}\right)  ,\ \ j\geq1.
\]
Thus, the UL-type $RG$-factorization is given by
\begin{equation}
I-\mathbb{P}=\left(  I-R_{U}\right)  \left(  I-\Phi_{D}\right)  \left(
I-G_{L}\right)  , \label{D-4}%
\end{equation}
where%
\[
R_{U}=\left(
\begin{array}
[c]{ccccc}%
0 & R &  &  & \\
& 0 & R &  & \\
&  & 0 & R & \\
&  &  & \ddots & \ddots
\end{array}
\right)  ,
\]%
\[
\Phi_{D}=\text{diag}\left(  \widehat{\Psi}_{0},\widehat{\Psi},\widehat{\Psi
},\widehat{\Psi},\cdots\right)
\]
and%
\[
G_{L}=\left(
\begin{array}
[c]{ccccc}%
0 &  &  &  & \\
G_{1,0} & 0 &  &  & \\
G_{2,0} & G_{1} & 0 &  & \\
G_{3,0} & G_{2} & G_{1} & 0 & \\
\vdots & \vdots & \vdots & \vdots & \ddots
\end{array}
\right)  .
\]
It follows from (\ref{D-3}) and (\ref{D-4}) that%
\begin{align*}
\pi &  =\left(  x_{0}B_{0},0,0,0,\ldots\right)  \left(  I-\mathbb{P}\right)
_{\min}^{-1}\\
&  =\left(  x_{0}B_{0},0,0,0,\ldots\right)  \left(  I-G_{L}\right)
^{-1}\left(  I-\Phi_{D}\right)  ^{-1}\left(  I-R_{U}\right)  ^{-1}\\
&  =\left(  x_{0}B_{0}\left(  I-\widehat{\Psi}_{0}\right)  ^{-1}%
,0,0,0,\ldots\right)  \left(  I-R_{U}\right)  ^{-1},
\end{align*}
where $\left(  I-\mathbb{P}\right)  _{\min}^{-1}=\sum_{k=0}^{\infty}%
\mathbb{P}^{k}$. Note that%
\[
\left(  I-R_{U}\right)  ^{-1}=\left(
\begin{array}
[c]{ccccc}%
I & R & R^{2} & R^{3} & \cdots\\
& I & R & R^{2} & \cdots\\
&  & I & R & \cdots\\
&  &  & I & \cdots\\
&  &  &  & \ddots
\end{array}
\right)  ,
\]
we obtain%
\begin{equation}
\pi_{k}=x_{0}B_{0}\left(  I-\widehat{\Psi}_{0}\right)  ^{-1}R^{k-1},\text{
\ }k\geq1. \label{D-5}%
\end{equation}

Comparing (\ref{D-5}) with (\ref{D-2}), we obtain%
\[
x_{0}B_{0}\left(  I-\widehat{\Psi}_{0}\right)  ^{-1}=x_{0}R_{1}\left(
I-R\right)  ^{-1}.
\]

\subsection{Markov chains of M/G/1 type}

We consider a discrete-time Markov chain $P$ of M/G/1 type whose transition
matrix is given by%
\begin{equation}
P=\left(
\begin{array}
[c]{ccccc}%
B_{1} & B_{2} & B_{3} & B_{4} & \cdots\\
B_{0} & A_{1} & A_{2} & A_{3} & \cdots\\
& A_{0} & A_{1} & A_{2} & \cdots\\
&  & A_{0} & A_{1} & \cdots\\
&  &  & \ddots & \ddots
\end{array}
\right)  . \label{D-6}%
\end{equation}
where the sizes of the two matrices $B_{1}$ and $A_{1}$ are $m_{0}$ and $m$,
respectively, while the sizes of other matrices can be determined accordingly.
We assume that this Markov chain is irreducible, aperiodic and positive
recurrent. Let the matrix $G$ be the minimal nonnegative solution to the
nonlinear matrix equation $G=\sum_{k=0}^{\infty}A_{k}G^{k}$.

In what follows we provide two methods to derive the tail probabilities in the
stationary regime.

\textbf{(a) The matrix-iterative solution}

Using Chapter 2 of Li \cite{Li:2010}, the $U$-measure is given by%
\[
\Psi_{0}=B_{1}+\sum_{k=2}^{\infty}B_{k}G^{k-2}G_{1}%
\]
and for $k\geq1$%
\[
\Psi=\Psi_{k}=\sum_{k=1}^{\infty}A_{k}G^{k-1};
\]
and the $R$-measure%
\[
R_{0,j}=\left(  \sum_{k=j+1}^{\infty}B_{k}G^{k-1}\right)  \left(
I-\Psi\right)  ^{-1},\ \ j\geq1,
\]
and%
\[
R_{j}=\left(  \sum_{k=j+1}^{\infty}A_{k}G^{k-1}\right)  \left(  I-\Psi\right)
^{-1},\ \ j\geq1.
\]
The stationary probability vector $x=\left(  x_{0},x_{1},x_{2},\ldots\right)
$ is given by%
\[
\left\{
\begin{array}
[c]{l}%
x_{0}=\tau y_{0},\\
x_{k}=x_{0}R_{0,k}+\sum\limits_{i=1}^{k-1}x_{i}R_{k-i},\ \ k\geq1,
\end{array}
\right.
\]
where $y_{0}$ is the stationary probability vector of the censored Markov
chain $\Psi_{0}$ to level $0$\ and the scalar $\tau$ is determined by
$\sum_{k=0}^{\infty}x_{k}e=1$ uniquely. Thus, we obtain%
\begin{equation}
\pi_{k}=\sum_{j=k}^{\infty}x_{j}=x_{0}\sum\limits_{j=k}^{\infty}R_{0,j}%
+\sum\limits_{i=1}^{\infty}x_{i}\sum\limits_{j=k}^{\infty}R_{j-i},\ \ k\geq1.
\label{D-7}%
\end{equation}

\textbf{(b) The UL-type }$\mathbf{RG}$\textbf{-factorization}

Note that $x=xP$, for $k\geq1$ we obtain%
\[
\pi_{k}=x_{0}\sum\limits_{j=k+1}^{\infty}B_{j}+\sum\limits_{i=2}^{k+1}\pi
_{i}A_{k+1-i}+\pi_{1}\sum\limits_{j=k}^{\infty}A_{j},
\]
this gives%
\begin{equation}
\pi=\pi\mathbf{P}+\left(  x_{0}\sum\limits_{j=2}^{\infty}B_{j},x_{0}%
\sum\limits_{j=3}^{\infty}B_{j},x_{0}\sum\limits_{j=4}^{\infty}B_{j}%
,\ldots\right)  , \label{D-8}%
\end{equation}
where%
\[
\mathbf{P=}\left(
\begin{array}
[c]{cccc}%
\sum\limits_{j=1}^{\infty}A_{j} & \sum\limits_{j=2}^{\infty}A_{j} &
\sum\limits_{j=3}^{\infty}A_{j} & \cdots\\
A_{0} & A_{1} & A_{2} & \cdots\\
& A_{0} & A_{1} & \cdots\\
&  & A_{0} & \cdots\\
&  &  & \ddots
\end{array}
\right)  .
\]
Using Chapter 2 of Li \cite{Li:2010}, the $U$-measure is given by%
\[
\widehat{\Psi}_{0}=\sum\limits_{j=1}^{\infty}A_{j}+\sum_{k=2}^{\infty}%
\sum\limits_{j=k}^{\infty}A_{j}G^{k-2}G_{1}%
\]
and for $k\geq1$%
\[
\widehat{\Psi}=\widehat{\Psi}_{k}=\sum_{i=1}^{\infty}A_{i}G^{i-1};
\]
and the $R$-measure%
\[
R_{0,j}=\left(  \sum_{k=j+1}^{\infty}B_{k}G^{k-1}\right)  \left(
I-\widehat{\Psi}\right)  ^{-1},\ \ j\geq1,
\]
and%
\[
R_{j}=\left(  \sum_{k=j+1}^{\infty}A_{k}G^{k-1}\right)  \left(  I-\widehat
{\Psi}\right)  ^{-1},\ \ j\geq1.
\]
Thus, the UL-type $RG$-factorization is given by
\begin{equation}
I-\mathbf{P}=\left(  I-R_{U}\right)  \left(  I-\Psi_{D}\right)  \left(
I-G_{L}\right)  , \label{D-9}%
\end{equation}
where%
\[
R_{U}=\left(
\begin{array}
[c]{ccccc}%
0 & R_{0,1} & R_{0,2} & R_{0,3} & \cdots\\
& 0 & R_{1} & R_{2} & \cdots\\
&  & 0 & R_{1} & \cdots\\
&  &  & 0 & \cdots\\
&  &  &  & \ddots
\end{array}
\right)  ,
\]%
\[
\Psi_{D}=\text{diag}\left(  \widehat{\Psi}_{0},\widehat{\Psi},\widehat{\Psi
},\widehat{\Psi},\ldots\right)
\]
and%
\[
G_{L}=\left(
\begin{array}
[c]{ccccc}%
0 &  &  &  & \\
G_{1} & 0 &  &  & \\
& G & 0 &  & \\
&  & G & 0 & \\
&  &  & \ddots & \ddots
\end{array}
\right)  .
\]
It follows from (\ref{D-8}) and (\ref{D-9}) that%
\begin{align*}
\pi &  =\left(  x_{0}\sum\limits_{j=2}^{\infty}B_{j},x_{0}\sum\limits_{j=3}%
^{\infty}B_{j},x_{0}\sum\limits_{j=4}^{\infty}B_{j},\ldots\right)  \left(
I-\mathbf{P}\right)  _{\min}^{-1}\\
&  =\left(  x_{0}\sum\limits_{j=2}^{\infty}B_{j},x_{0}\sum\limits_{j=3}%
^{\infty}B_{j},x_{0}\sum\limits_{j=4}^{\infty}B_{j},\ldots\right)  \left(
I-G_{L}\right)  ^{-1}\left(  I-\Psi_{D}\right)  ^{-1}\left(  I-R_{U}\right)
^{-1}.
\end{align*}

\section{Some Queueing Examples}

In this section, we consider four queueing examples which indicate how to use
our above results. We first provide a detailed discussion for the M/M/1
retrial queue with exponentially distributed retrial times. Then we simply
analyze other three queueing examples: The M(n)/M(n)/1 queue, the M/M/1 queue
with server multiple vacations, and the M/M/1 queue with repairable server.

\subsection{The M/M/1 retrial queue}

We consider an M/M/1 retrial queue with exponentially distributed retrial
times, where the arrival, service and retrial rates are $\lambda$, $\mu$ and
$\theta$, respectively. We denote by $N\left(  t\right)  $ and $C\left(
t\right)  $ the number of customers in the orbit and the state of server at
time $t$, respectively, where $N\left(  t\right)  =0,1,2,\ldots$ and $C\left(
t\right)  =W$ for the busy server or $I$ for the idle server. For $k\geq0$, we
write%
\[
p_{W,k}\left(  t\right)  =P\left\{  C\left(  t\right)  =W,N\left(  t\right)
=k\right\}
\]
and%
\[
p_{I,k}\left(  t\right)  =P\left\{  C\left(  t\right)  =I,N\left(  t\right)
=k\right\}  .
\]
Hence, we obtain%
\begin{align}
\frac{\text{d}}{\text{d}t}p_{W,0}\left(  t\right)   &  =-\left(  \lambda
+\mu\right)  p_{W,0}\left(  t\right)  +\lambda p_{I,0}\left(  t\right)
+\theta p_{I,1}\left(  t\right)  ,\nonumber\\
\frac{\text{d}}{\text{d}t}p_{W,k}\left(  t\right)   &  =-\left(  \lambda
+\mu\right)  p_{W,k}\left(  t\right)  +\lambda p_{I,k}\left(  t\right)
+\left(  k+1\right)  \theta p_{I,k+1}\left(  t\right)  +\lambda p_{W,k-1}%
\left(  t\right)  ,k\geq1, \label{Equation1}%
\end{align}%
\begin{align}
\frac{\text{d}}{\text{d}t}p_{I,0}\left(  t\right)   &  =\mu p_{W,0}\left(
t\right)  -\lambda p_{I,0}\left(  t\right)  ,\nonumber\\
\frac{\text{d}}{\text{d}t}p_{I,k}\left(  t\right)   &  =\mu p_{W,k}\left(
t\right)  -\lambda p_{I,k}\left(  t\right)  -k\theta p_{I,k}\left(  t\right)
,\text{ \ }k\geq1. \label{Equation2}%
\end{align}
Let $\rho=\lambda/\mu<1$. Then the M/M/1 retrial queue is stable. In this
case, we write that for $k\geq0$%
\[
x_{W,k}=\lim_{t\rightarrow+\infty}p_{W,k}\left(  t\right)  ,\text{ \ }%
x_{I,k}=\lim_{t\rightarrow+\infty}p_{I,k}\left(  t\right)  ,
\]
and%
\[
\pi_{W,k}=\sum_{j=k}^{\infty}x_{W,j},\text{ \ \ }\pi_{I,k}=\sum_{j=k}^{\infty
}x_{I,j}.
\]
Then it follows from (\ref{Equation1}) and (\ref{Equation2}) that%
\begin{equation}
\mu\pi_{W,0}-\lambda\pi_{I,0}-\theta\sum_{j=1}^{\infty}\pi_{I,j}=0,
\label{FP5}%
\end{equation}%
\begin{equation}
\pi_{W,0}+\pi_{I,0}=1, \label{FP6}%
\end{equation}
and for $k\geq1$%
\begin{equation}
\lambda\left(  \pi_{W,k-1}-\pi_{W,k}\right)  -\mu\pi_{W,k}+\lambda\pi
_{I,k}+\theta\left[  \left(  k+1\right)  \pi_{I,k+1}+\sum_{j=k+2}^{\infty}%
\pi_{I,j}\right]  =0 \label{FP7}%
\end{equation}
and%
\begin{equation}
\mu\pi_{W,k}-\lambda\pi_{I,k}-\theta\left[  k\pi_{I,k}+\sum_{j=k+1}^{\infty
}\pi_{I,j}\right]  =0. \label{FP8}%
\end{equation}

Let%
\[
\Pi=\left(  \Pi_{1},\Pi_{2},\Pi_{3},\ldots\right)  ,
\]%
\[
\Pi_{k}=\left(  \pi_{W,k},\pi_{I,k}\right)  ,\text{ \ }k\geq1;
\]%
\begin{equation}
Q=\left(
\begin{array}
[c]{ccccccc}%
A_{1} & C &  &  &  &  & \\
B_{2} & A_{2} & C &  &  &  & \\
D & B_{3} & A_{3} & C &  &  & \\
D & D & B_{4} & A_{4} & C &  & \\
D & D & D & B_{5} & A_{5} & C & \\
\vdots & \vdots & \vdots & \vdots & \vdots & \vdots & \ddots
\end{array}
\right)  , \label{FP14}%
\end{equation}
and for \ $k\geq1$%
\[
A_{k}=\left(
\begin{array}
[c]{cc}%
-\left(  \lambda+\mu\right)  & \mu\\
\lambda & -\left(  \lambda+k\theta\right)
\end{array}
\right)  ,\text{ \ }B_{k+1}=\left(
\begin{array}
[c]{cc}%
0 & 0\\
\left(  k+1\right)  \theta & -\theta
\end{array}
\right)  ,
\]%
\[
C=\left(
\begin{array}
[c]{cc}%
\lambda & 0\\
0 & 0
\end{array}
\right)  ,\text{ \ }D=\left(
\begin{array}
[c]{cc}%
0 & 0\\
\theta & -\theta
\end{array}
\right)  .
\]
Note that $\pi_{W,0}=\rho$ and $\pi_{I,0}=1-\rho$, it follows from (\ref{FP7})
and (\ref{FP8}) that%
\begin{equation}
\Pi Q=\left(  -\lambda\rho,0,0,\ldots\right)  . \label{FP15}%
\end{equation}

To solve Equation (\ref{FP15}), we need to construct a UL-type $RG$%
-factorization of the matrix $Q$ in which the computational steps are similar
to that in Subsection 2.2.3 of Li \cite{Li:2010}. Here, we provide a sketch of
the computation as follows. Let%
\[
W_{k}=\left(
\begin{array}
[c]{cccccc}%
A_{k} & C &  &  &  & \\
B_{k+1} & A_{k+1} & C &  &  & \\
D & B_{k+2} & A_{k+2} & C &  & \\
D & D & B_{k+3} & A_{k+3} & C & \\
\vdots & \vdots & \vdots & \vdots & \vdots & \ddots
\end{array}
\right)  .
\]
We denote by $\left(  \widehat{W}_{1,1}^{\left(  k\right)  },\widehat{W}%
_{1,2}^{\left(  k\right)  },\widehat{W}_{1,3}^{\left(  k\right)  }%
,\ldots\right)  $ the first block-row of the matrix $\left(  -W_{k}\right)
_{\min}^{-1}$. Thus for $k\geq1$ and $j\geq2$,%
\begin{align*}
R_{k,k+1}  &  =C\widehat{W}_{1,1}^{\left(  k+1\right)  }\overset{\text{def}%
}{=}R_{k},\\
R_{k,k+j}  &  =0;
\end{align*}
and for $i\geq2$ and $1\leq j\leq i-2$%
\begin{align*}
G_{i,i-1}  &  =\widehat{W}_{1,1}^{\left(  k\right)  }B_{i}+\left[  \sum
_{l=2}^{\infty}\widehat{W}_{1,l}^{\left(  k\right)  }\right]  D\\
G_{i,j}  &  =\left[  \sum_{l=1}^{\infty}\widehat{W}_{1,l}^{\left(  k\right)
}\right]  D\overset{\text{def}}{=}G_{i}.
\end{align*}
In what follows we provide some further interpretation on the $R$- and
$G$-measures. Let the $R$-measure $\left\{  R_{k}:k\geq1\right\}  $ be the
minimal nonnegative solution to the following system of nonlinear equations%
\[
C+R_{k}A_{k}+R_{k}R_{k+1}B_{k+2}+R_{k}R_{k+1}\left(  \sum_{l=2}^{\infty
}R_{k+2}R_{k+3}\cdots R_{k+l}\right)  D=0,\text{ \ }k\geq1.
\]
Once the $R$-measure $\left\{  R_{k}:k\geq1\right\}  $ is determined, we have%
\[
\Psi_{k}=A_{k}+R_{k}B_{k}+R_{k}\left(  \sum_{l=1}^{\infty}R_{k+1}R_{k+2}\cdots
R_{k+l}\right)  D,
\]%
\[
\widehat{W}_{1,1}^{\left(  k\right)  }=\left(  -\Psi_{k}\right)  ^{-1},
\]%
\[
\widehat{W}_{1,j}^{\left(  k\right)  }=\widehat{W}_{1,1}^{\left(  k\right)
}R_{k}R_{k+1}\cdots R_{k+j-2}=\left(  -\Psi_{k}\right)  ^{-1}R_{k}%
R_{k+1}\cdots R_{k+j-2},\text{ \ }j\geq2;
\]
and for $i\geq2$ and $1\leq j\leq i-2$%
\begin{align*}
G_{i,i-1}  &  =\left(  -\Psi_{k}\right)  ^{-1}\left[  B_{i}+\left(  \sum
_{l=1}^{\infty}R_{i}R_{i+1}\cdots R_{i+l-1}\right)  D\right]  ,\\
G_{i,j}  &  =\left(  -\Psi_{i}\right)  ^{-1}\left(  I+\sum_{l=1}^{\infty}%
R_{i}R_{i+1}\cdots R_{i+l-1}\right)  D\overset{\text{def}}{=}G_{i}.
\end{align*}
Thus the UL-type $RG$-factorization is given by%
\begin{equation}
Q=\left(  I-R_{U}\right)  U_{D}\left(  I-G_{L}\right)  , \label{FP17}%
\end{equation}
where%
\[
R_{U}=\left(
\begin{array}
[c]{cccccc}%
0 & R_{1} &  &  &  & \\
& 0 & R_{2} &  &  & \\
&  & 0 & R_{3} &  & \\
&  &  & 0 & R_{4} & \\
&  &  &  & \ddots & \ddots
\end{array}
\right)  ,
\]%
\[
U_{D}=\text{diag}\left(  \Psi_{1},\Psi_{2},\Psi_{3},\Psi_{4},\ldots\right)
\]
and%
\[
G_{L}=\left(
\begin{array}
[c]{ccccc}%
0 &  &  &  & \\
G_{2,1} & 0 &  &  & \\
G_{3} & G_{3,2} & 0 &  & \\
G_{4} & G_{4} & G_{4,3} & 0 & \\
\vdots & \vdots & \vdots & \vdots & \ddots
\end{array}
\right)  .
\]
It follows from (\ref{FP15}) and (\ref{FP17}) that%
\begin{align*}
\Pi &  =\left(  -\lambda\rho,0,0,\ldots\right)  \left(  I-G_{L}\right)
^{-1}U_{D}^{-1}\left(  I-R_{U}\right)  ^{-1}\\
&  =\left(  \left(  -\lambda\rho,0\right)  \Psi_{1}^{-1},0,0,\ldots\right)
\left(  I-R_{U}\right)  ^{-1}.
\end{align*}
Note that%
\[
\left(  I-R_{U}\right)  ^{-1}=\ \left(
\begin{array}
[c]{ccccc}%
I & R_{1} & R_{1}R_{2} & R_{1}R_{2}R_{3} & \cdots\\
& I & R_{1} & R_{1}R_{2} & \cdots\\
&  & I & R_{1} & \cdots\\
&  &  & I & \cdots\\
&  &  &  & \ddots
\end{array}
\right)  ,
\]
we obtain%
\begin{align*}
\pi_{1}  &  =\left(  \lambda\rho,0\right)  \left(  -\Psi_{1}^{-1}\right)  ,\\
\pi_{k}  &  =\left(  \lambda\rho,0\right)  \left(  -\Psi_{1}^{-1}\right)
R_{1}R_{2}\cdots R_{k-1},\text{ \ }k\geq2.
\end{align*}

\subsection{The M(n)/M(n)/1 queue}

We consider an M(n)/M(n)/1 queue whose arrival and service rates depend on the
number of customers in this system, denoted as $\lambda_{n}$ and $\mu_{n}$,
respectively. We denote by $N\left(  t\right)  $ the number of customers in
this system at time $t$. Then $N\left(  t\right)  \in\left\{  0,1,2,\ldots
\right\}  $. For $k\geq0$, we write%
\[
Q_{k}\left(  t\right)  =P\left\{  N\left(  t\right)  \geq k\right\}
\]
and when the M(n)/M(n)/1 queue is stable,
\[
\pi_{k}=\lim_{t\rightarrow+\infty}Q_{k}\left(  t\right)  .
\]
Then we obtain that for $k\geq1$%
\[
\lambda_{k-1}\left(  \pi_{k-1}-\pi_{k}\right)  =\mu_{k}\left(  \pi_{k}%
-\pi_{k+1}\right)
\]
with the boundary condition $\pi_{0}=1$. Let $\pi_{1}=g\in\left(  0,1\right)
$ and $\rho_{k-1}=\lambda_{k-1}/\mu_{k}$ for $k\geq1$. Then%
\[
\pi_{1}-\pi_{2}=\rho_{0}\left(  \pi_{0}-\pi_{1}\right)  =\rho_{0}\left(
1-g\right)
\]
and for $k\geq2$%
\[
\pi_{k}-\pi_{k+1}=\rho_{k-1}\left(  \pi_{k-1}-\pi_{k}\right)  =\rho_{k-1}%
\rho_{k-2}\cdots\rho_{1}\rho_{0}\left(  1-g\right)  .
\]
We obtain%
\[
g=\left(  \rho_{0}+\rho_{1}\rho_{0}+\rho_{2}\rho_{1}\rho_{0}+\rho_{3}\rho
_{2}\rho_{1}\rho_{0}+\cdots\right)  \left(  1-g\right)  ,
\]
from which follows%
\[
g=\frac{\rho_{0}+\rho_{1}\rho_{0}+\rho_{2}\rho_{1}\rho_{0}+\rho_{3}\rho
_{2}\rho_{1}\rho_{0}+\cdots}{1+\left(  \rho_{0}+\rho_{1}\rho_{0}+\rho_{2}%
\rho_{1}\rho_{0}+\rho_{3}\rho_{2}\rho_{1}\rho_{0}+\cdots\right)  }.
\]
Thus for $k\geq1$%
\[
\pi_{k}=\frac{\rho_{k-1}\rho_{k-2}\cdots\rho_{1}\rho_{0}+\rho_{k}\rho
_{k-1}\rho_{k-2}\cdots\rho_{1}\rho_{0}+\rho_{k+1}\rho_{k}\rho_{k-1}\rho
_{k-2}\cdots\rho_{1}\rho_{0}+\cdots}{1+\left(  \rho_{0}+\rho_{1}\rho_{0}%
+\rho_{2}\rho_{1}\rho_{0}+\rho_{3}\rho_{2}\rho_{1}\rho_{0}+\cdots\right)  }.
\]

It is interesting to extend the above result to more general models such as
the MAP(n)/M/1 queue and the M/PH(n)/1 queue. The more general queues can be
analyzed by the level-dependent QBD processes, see Section 3.

\subsection{The M/M/1 queue with server multiple vacations}

We consider an M/M/1 queue with server multiple vacations, where the arrival,
service and vacation rates are $\lambda$, $\mu=1$ and $\theta$. The vacation
process is based on the multiple vacation policy: When there is no customer in
the system, the server immediately proceeds on vacation and keeps taking
vacations until it finds at least one customer waiting in the server or its
buffer at the vacation completion instant. The arrival, service and vacation
processes are independent of each other.

Let $N\left(  t\right)  $ be the number of customers in the queueing system at
time $t$, and%
\[
\xi\left(  t\right)  =\left\{
\begin{array}
[c]{ll}%
V, & \text{if the server is taking a vacation at time }t,\\
W, & \text{if the server is working at time }t.
\end{array}
\right.
\]
Then $\left\{  \left(  \xi\left(  t\right)  ,N\left(  t\right)  \right)
:t\geq0\right\}  $ is a Markov chain on a state space $E=\{\left(  V,k\right)
,\left(  W,l\right)  :$\ $k\geq0,l\geq1\}$. We write%
\[
Q_{V,l}\left(  t\right)  =P\left\{  \xi\left(  t\right)  =V,N\left(  t\right)
\geq l\right\}  ,\text{ \ }l\geq0,
\]
and%
\[
Q_{W,k}\left(  t\right)  =P\left\{  \xi\left(  t\right)  =W,N\left(  t\right)
\geq k\right\}  ,\text{ \ }k\geq1.
\]
If $0<\lambda<\mu=1$, then this queue is stable. We set%
\[
\pi_{V,k}=\lim_{t\rightarrow+\infty}Q_{V,k}\left(  t\right)  ,\text{
\ \ }k\geq0,
\]%
\[
\pi_{W,l}=\lim_{t\rightarrow+\infty}Q_{W,l}\left(  t\right)  ,\text{
\ \ }l\geq1.
\]
Then we obtain%
\begin{equation}
\left(  \pi_{W,1}-\pi_{W,2}\right)  -\theta\pi_{V,1}=0,\label{EquT10}%
\end{equation}%
\begin{equation}
\lambda\left(  \pi_{V,k-1}-\pi_{V,k}\right)  -\theta\pi_{V,k}=0,\text{
\ }k\geq1,\label{EquT11}%
\end{equation}%
\begin{equation}
\lambda\left(  \pi_{W,l-1}-\pi_{W,l}\right)  -\left(  \pi_{W,l}-\pi
_{W,l+1}\right)  +\theta\pi_{V,l}=0,\text{ \ \ }l\geq2.\label{EquT12}%
\end{equation}
Note that $\pi_{V,0}=1-\lambda$ and $\pi_{W,1}=\lambda$, thus it follows from
(\ref{EquT11}) that%
\[
\pi_{V,k}=\left(  \frac{\lambda}{\lambda+\theta}\right)  ^{k}\left(
1-\lambda\right)  ,\text{ \ }k\geq0,
\]
and from (\ref{EquT10}) that%
\[
\pi_{W,2}=\lambda-\frac{\lambda\theta}{\lambda+\theta}\left(  1-\lambda
\right)  .
\]
Using $\pi_{W,1}=\lambda$ and $\pi_{W,2}=\lambda-\lambda\theta\left(
1-\lambda\right)  /\left(  \lambda+\theta\right)  $, it follows from
(\ref{EquT12}) that for $k\geq3$%
\begin{equation}
\pi_{W,k}=\pi_{W,k-1}-\lambda\left(  \pi_{W,k-2}-\pi_{W,k-1}\right)
-\theta\pi_{V,k-1},\label{EquT13}%
\end{equation}
which can be computed iteratively.

Let%
\[
Q=\left(
\begin{array}
[c]{ccccc}%
-\left(  1+\lambda\right)  & \lambda &  &  & \\
1 & -\left(  1+\lambda\right)  & \lambda &  & \\
& 1 & -\left(  1+\lambda\right)  & \lambda & \\
&  & \ddots & \ddots & \ddots
\end{array}
\right)  .
\]
Then using (\ref{EquT13}) we obtain%
\begin{equation}
\left(  \pi_{W,2},\pi_{W,3},\pi_{W,4},\ldots\right)  Q=-\left(  \lambda
\pi_{W,1}+\theta\pi_{V,2},\theta\pi_{V,3},\theta\pi_{V,4},\ldots\right)  .
\label{EquT14}%
\end{equation}
Let%
\[
R=\lambda,\text{ \ \ }G=1.
\]
Then%
\[
U_{k}=-\left(  1+\lambda\right)  +R=-1,\text{ \ \ }k\geq0,
\]%
\[
U_{D}=\text{diag}\left(  -1,-1,-1,-1,\ldots\right)  ,
\]%
\[
R_{U}=\left(
\begin{array}
[c]{ccccc}%
0 & \lambda &  &  & \\
& 0 & \lambda &  & \\
&  & 0 & \lambda & \\
&  &  & \ddots & \ddots
\end{array}
\right)
\]
and%
\[
G_{L}=\left(
\begin{array}
[c]{cccc}%
0 &  &  & \\
1 & 0 &  & \\
& 1 & 0 & \\
&  & \ddots & \ddots
\end{array}
\right)  .
\]
Thus we obtain%
\begin{align*}
\left(  \pi_{W,2},\pi_{W,3},\pi_{W,4},\ldots\right)   &  =-\left(  \lambda
\pi_{W,1}+\theta\pi_{V,2},\theta\pi_{V,3},\theta\pi_{V,4},\ldots\right)
Q_{\max}^{-1}\\
&  =-\left(  \lambda\pi_{W,1}+\theta\pi_{V,2},\theta\pi_{V,3},\theta\pi
_{V,4},\ldots\right)  \left(  I-G_{L}\right)  ^{-1}U_{D}^{-1}\left(
I-R_{U}\right)  ^{-1}.
\end{align*}
Note that%
\[
\left(  I-G_{L}\right)  ^{-1}=\left(
\begin{array}
[c]{ccccc}%
1 &  &  &  & \\
1 & 1 &  &  & \\
1 & 1 & 1 &  & \\
1 & 1 & 1 & 1 & \\
\vdots & \vdots & \vdots & \vdots & \ddots
\end{array}
\right)
\]
and%
\[
\left(  I-R_{U}\right)  ^{-1}=\left(
\begin{array}
[c]{ccccc}%
1 & \lambda & \lambda^{2} & \lambda^{3} & \cdots\\
& 1 & \lambda & \lambda^{2} & \cdots\\
&  & 1 & \lambda & \cdots\\
&  &  & 1 & \cdots\\
&  &  &  & \ddots
\end{array}
\right)  ,
\]
for $k\geq2$ we can obtain%
\begin{align*}
\pi_{W,k}  &  =\lambda^{k-1}\pi_{W,1}+\theta\sum_{i=0}^{k-2}\lambda^{i}%
\sum_{l=k-i}^{\infty}\pi_{V,l}\\
&  =\lambda^{k}+\left(  1-\lambda\right)  \frac{\lambda^{2}}{\theta}\left[
1-\left(  \frac{\lambda}{\lambda+\theta}\right)  ^{k-1}\right]  .
\end{align*}

\subsection{The M/M/1 queue with repairable server}

We consider an M/M/1 queue with repairable server, where the arrival and
service rates are $\lambda$ and $\mu$, respectively. The life time of the
server is exponential with failure rate $\alpha$. Once the server failed, it
immediately is repaired, and the repair time is exponential with repair rate
$\beta$. The repaired server is the same as the new one. We assume that all
the random variables defined above are independent of each other.

For this M/M/1 repairable queue, we denote by $N\left(  t\right)  $ and
$C\left(  t\right)  $ the number of customers in this queueing system and the
state of the server at time $t\geq0$, respectively, where $N\left(  t\right)
=0,1,2,\ldots$, and $C\left(  t\right)  =W$ for server working or $R$ for
server repair. It is easy to see that $\left\{  \left(  N\left(  t\right)
,C\left(  t\right)  \right)  :t\geq0\right\}  $ is a Markov chain. For
$k\geq0$ and $l\geq1$, we write%
\[
Q_{W,k}\left(  t\right)  =P\left\{  C\left(  t\right)  =W,N\left(  t\right)
\geq k\right\}
\]
and%
\[
Q_{R,l}\left(  t\right)  =P\left\{  C\left(  t\right)  =R,N\left(  t\right)
\geq l\right\}  .
\]
If $\rho=\frac{\lambda}{\mu}\left(  1+\frac{\alpha}{\beta}\right)  <1$, then
this queue is stable. Let%
\[
\pi_{W,k}=\lim_{t\rightarrow+\infty}Q_{W,k}\left(  t\right)  ,\text{ \ }%
k\geq0,
\]
and%
\[
\pi_{R,l}=\lim_{t\rightarrow+\infty}Q_{R,l}\left(  t\right)  ,\text{ \ }%
l\geq1.
\]
Then we obtain%
\begin{equation}
\pi_{W,0}+\pi_{R,1}=1\label{Equa8-1}%
\end{equation}%
\begin{equation}
-\alpha\pi_{W,1}+\beta\pi_{R,1}=0,\label{Equa9}%
\end{equation}
for $k\geq1$%
\begin{equation}
\lambda\left(  \pi_{W,k-1}-\pi_{W,k}\right)  -\mu\left(  \pi_{W,k}-\pi
_{W,k+1}\right)  -\alpha\pi_{W,k}+\beta\pi_{R,k}=0,\label{Equa10}%
\end{equation}
for $l\geq2$%
\begin{equation}
\lambda\left(  \pi_{R,l-1}-\pi_{R,l}\right)  +\alpha\pi_{W,l}-\beta\pi
_{R,l}=0.\label{Equa11}%
\end{equation}

It follows from (\ref{Equa10}) that%
\begin{equation}
\lambda\pi_{W,0}-\mu\pi_{W,1}-\alpha\sum_{k=1}^{\infty}\pi_{W,k}+\beta
\sum_{k=1}^{\infty}\pi_{R,k}=0 \label{Equa12}%
\end{equation}
and from (\ref{Equa11}) that%
\[
\lambda\pi_{R,1}+\alpha\sum_{k=2}^{\infty}\pi_{W,k}-\beta\sum_{k=2}^{\infty
}\pi_{R,k}=0,
\]
which, together with (\ref{Equa9}), leads to%
\begin{equation}
\lambda\pi_{R,1}+\alpha\sum_{k=1}^{\infty}\pi_{W,k}-\beta\sum_{k=1}^{\infty
}\pi_{R,k}=0. \label{Equa13}%
\end{equation}
Using (\ref{Equa12}) and (\ref{Equa13}), we obtain%
\begin{equation}
\lambda\pi_{W,0}-\mu\pi_{W,1}+\lambda\pi_{R,1}=0. \label{Equa14}%
\end{equation}
It follows from (\ref{Equa8-1}), (\ref{Equa9}) and (\ref{Equa14}) that%
\[
\pi_{W,0}=1-\frac{\lambda}{\mu}\frac{\alpha}{\beta},
\]%
\[
\pi_{W,1}=\frac{\lambda}{\mu}%
\]
and%
\[
\pi_{R,1}=\frac{\lambda}{\mu}\frac{\alpha}{\beta}.
\]
It follows from (\ref{Equa10}) and (\ref{Equa11}) that for $k\geq2$%
\begin{equation}
\pi_{W,k}=\frac{\lambda+\mu+\alpha}{\mu}\pi_{W,k-1}-\frac{\lambda}{\mu}%
\pi_{W,k-2}-\frac{\beta}{\mu}\pi_{R,k-1} \label{Equa15}%
\end{equation}
and%
\begin{equation}
\pi_{R,k}=\frac{\alpha}{\lambda+\beta}\pi_{W,k}+\frac{\lambda}{\lambda+\beta
}\pi_{R,k-1}. \label{Equa16}%
\end{equation}
Therefore, $\pi_{W,k}$ and $\pi_{R,k}$ for $k\geq2$ can be computed iteratively.

To provide explicit expressions for $\pi_{W,k}$ and $\pi_{R,k}$ with $k\geq2$,
we write%
\[
\Pi_{k}=\left(  \pi_{W,k},\pi_{R,k}\right)  ,\text{ \ }k\geq2,
\]%
\[
\Pi=\left(  \Pi_{2},\Pi_{3},\Pi_{4},\Pi_{5},\ldots\right)  ,
\]%
\[
A=\left(
\begin{array}
[c]{cc}%
-\left(  \lambda+\mu+\alpha\right)   & \alpha\\
\beta & -\left(  \lambda+\beta\right)
\end{array}
\right)  ,\text{ }B=\left(
\begin{array}
[c]{cc}%
\mu & 0\\
0 & 0
\end{array}
\right)  ,\text{ }C=\left(
\begin{array}
[c]{cc}%
\lambda & \\
& \lambda
\end{array}
\right)  ,
\]%
\[
\mathcal{Q}=\left(
\begin{array}
[c]{ccccc}%
A & C &  &  & \\
B & A & C &  & \\
& B & A & C & \\
&  & \ddots & \ddots & \ddots
\end{array}
\right)  .
\]
It follows from (\ref{Equa15}) and (\ref{Equa16}) that%
\begin{equation}
\Pi\mathcal{Q}=-\left(  \left(  \frac{\lambda^{2}}{\mu},\frac{\lambda^{2}}%
{\mu}\frac{\alpha}{\beta}\right)  ,0,0,0,\ldots\right)  .\label{Equa17}%
\end{equation}
Let $R$ and $G$\ be the minimal nonnegative solutions to the nonlinear
equations $C+RA+R^{2}B=0$ and $CG^{2}+AG+B=0$, respectively. It is easy to see
that the infinitesimal generator $\mathcal{Q}$ has the UL-type $RG$%
-factorization $\mathcal{Q}=\left(  I-R_{U}\right)  U_{D}\left(
I-G_{L}\right)  $, where%
\[
U_{D}=\text{diag}\left(  \Psi,\Psi,\Psi,\Psi,\ldots\right)  ,\text{ \ }%
\Psi=A+RB=A+CG,
\]%
\[
R_{U}=\left(
\begin{array}
[c]{ccccc}%
0 & R &  &  & \\
& 0 & R &  & \\
&  & 0 & R & \\
&  &  & \ddots & \ddots
\end{array}
\right)  ,\text{ \ }G_{L}=\left(
\begin{array}
[c]{ccccc}%
0 &  &  &  & \\
G & 0 &  &  & \\
& G & 0 &  & \\
&  & G & 0 & \\
&  &  & \ddots & \ddots
\end{array}
\right)  .
\]
Thus It follows from (\ref{Equa17}) that%
\begin{align*}
\Pi &  =\left(  \left(  \frac{\lambda^{2}}{\mu},\frac{\lambda^{2}}{\mu
}\frac{\alpha}{\beta}\right)  ,0,0,0,\ldots\right)  \left(  I-G_{L}\right)
^{-1}\left(  -U_{D}^{-1}\right)  \left(  I-R_{U}\right)  ^{-1}\\
&  =\left(  \left(  \frac{\lambda^{2}}{\mu},\frac{\lambda^{2}}{\mu
}\frac{\alpha}{\beta}\right)  \left(  -\Psi^{-1}\right)  ,0,0,0,\ldots\right)
\left(  I-R_{U}\right)  ^{-1}.
\end{align*}
This gives%
\begin{align*}
\Pi_{2} &  =\left(  \frac{\lambda^{2}}{\mu},\frac{\lambda^{2}}{\mu
}\frac{\alpha}{\beta}\right)  \left(  -\Psi^{-1}\right)  ,\\
\Pi_{k} &  =\left(  \frac{\lambda^{2}}{\mu},\frac{\lambda^{2}}{\mu
}\frac{\alpha}{\beta}\right)  \left(  -\Psi^{-1}\right)  R^{k-2},\text{
\ }k\geq3.
\end{align*}
In fact, the minimal nonnegative solution $R$ can be explicitly determined
from the nonlinear 2-order matrix equation $C+RA+R^{2}B=0$, here we omit the detail.

\section{Concluding remarks}

This paper discusses tail probabilities of queueing processes, such as, the
QBD processes and Markov chains of GI/M/1 type and of M/G/1 type, and provides
some efficient algorithms for computing the tail probabilities by means of the
matrix-geometric solution, the matrix-iterative solution, the matrix-product
solution and the two types of $RG$-factorizations. Also, we consider four
queueing examples: The M/M/1 retrial queue, the M(n)/M(n)/1 queue, the M/M/1
queue with server multiple vacations, and the M/M/1 queue with repairable
server, where the M/M/1 retrial queue is given a detailed discussion, while
for the other three queues, a sketch of the analysis is given. It is seen from
the four queueing examples that the method of this paper can be applied to
deal with more general queues including the MAP/PH/1 queue, the GI/PH/1 queue
and the BMAP/SM/1 queue.

The results given in this paper are very useful in the study of large scale
stochastic networks with resource management, such as, supermarket models and
work stealing models. Also, it will open a new avenue to helpfully analyze the
tail probabilities of many large scale stochastic networks when applying
differential equations and mean-field limits.

\section*{Acknowledgements}

The author thanks two reviewers for many valuable comments to sufficiently
improve the presentation of this paper. At the same time, the author
acknowledges that this research is partly supported by the National Natural
Science Foundation of China (No. 71271187) and the Hebei Natural Science
Foundation of China (No. A2012203125).

\end{document}